\documentclass[11pt,intlimits, draft,oneside]{amsart}

\usepackage[latin2]{inputenc}
\usepackage{latexsym}
\usepackage{amssymb}
\usepackage{amsthm}

\usepackage{amsmath}

\usepackage{enumerate}

\input{cyracc.def}
\newfont{\cyrfnt}{wncyi10 at 11pt}

\newtheorem{thm}{Theorem}[section]

\theoremstyle{definition}

\newtheorem{remark}[thm]{Remark}
\newtheorem{example}[thm]{Example}
\newtheorem{question}[thm]{Question}

\newcommand{\R}{\mathbb{R}}

\newcommand{\N}{\mathbb{N}}

\renewcommand{\Re}{\mathrm{Re}\,}

\begin{document}

\title[Weak analogue of the Trotter-Kato theorem]
{On the weak analogue of the Trotter-Kato theorem}

\author{Tanja Eisner}
\address{Tanja Eisner \newline Mathematisches Institut, Universit\"{a}t T\"{u}bingen\newline Auf der Morgenstelle 10, D-70176, T\"{u}bingen, Germany}
\email{talo@fa.uni-tuebingen.de}

\renewcommand{\thefootnote}{}
\author{Andr\'{a}s Ser\'{e}ny}
\footnote{The second author was supported by the DAAD-PPP-Hungary Grant, project number D/05/01422.}
\address{Andr\'{a}s Ser\'{e}ny \newline Department of Mathematics and its Applications, Central European University,\newline N\'{a}dor utca  9, H-1051 Budapest, Hungary }
\email{sandris@elte.hu}


\begin{abstract}
In the Trotter--Kato approximation theorem for $C_0$-semigroups on
Banach spaces, we replace the strong by the weak operator topology and
discuss the validity of the relevant implications. 
\end{abstract}
\maketitle



\section{Introduction}

Trotter--Kato approximation theorems relate the convergence of a
sequence of $C_0$-semigroups $(T_n(t))_{t\geq 0}$ on a Banach space 
to the convergence of the corresponding generators $A_n$ and
their resolvents $R(\lambda, A_n)$, in each case for the strong
operator topology. They 
are still an interesting topic in semigroup theory (see, e.g., Bobrowski \cite{bobrowski:2007})
and are used, e.g., for partial and stochastic (partial) differential
equations and in numerical analysis (see, e.g., Ito and Kappel \cite{ito/kappel:1998}).

In this note we consider the following question: Do these results
still hold for convergence with respect to the weak operator
topology?
%

First positive results have been proved by Yosida \cite[Section
IX.12]{yosida:1974} for semigroups on locally convex spaces assuming
equicontinuity. Later, using the concept of bi-continuous semigroups
introduced by K\"uhnemund \cite{kuehnemund:2003}, Albanese, K\"uhnemund
\cite{albanese/kuehnemund:2002} and Albanese, Mangino
\cite{albanese/mangino:2004} extended the Trotter--Kato theorems to
this situation and applied it to Markov semigroups. They
also needed equicontinuity, an assumption far too strong for
the weak operator topology. 

In this paper we 
discuss this situation and prove or disprove the
implications in the following Trotter--Kato theorem quoted from
Engel, Nagel \cite[Theorem III.4.8]{engel/nagel:2000} if we replace
the strong by the weak operator topology.
%
%
\begin{thm}\label{thm:trotter-kato-1}\emph{(First Trotter--Kato theorem)}
Let $T(\cdot)$ and $T_n(\cdot)$, $n\in\N$, be $C_0$-semigroups on a Banach space $X$ with generators $A$ and $A_n$, $n\in\N$, respectively. Assume that for some constants $M\geq 1$, $\omega\in\R$ the estimates $\|T(t)\|$,  $\|T_n(t)\|\leq Me^{\omega t}$ hold for all $t\geq 0$ and $n\in\N$. Let further $D$ be a core for $A$. Consider the following assertions.
\begin{enumerate}[(a)]
\item $D\subset D(A_n)$ for every $n\in\N$ and $\lim_{n\to\infty}A_n x=Ax$ for all $x\in D$.
\item For every $x\in D$ there exists $x_n\in D(A_n)$ such that $\lim_{n\to\infty}x_n=x$ and $\lim_{n\to\infty}A_n x_n=Ax$.
\item $\lim_{n\to\infty}R(\lambda, A_n)x=R(\lambda, A)x$ for every $x\in X$ and some $\lambda$ with $\Re\lambda>\omega$.
\item $\lim_{n\to\infty}T_n(t)x=T(t)x$ for every $x\in X$, uniformly in $t$ on compact intervals.
\end{enumerate}
Then $(a)\Rightarrow (b) \Leftrightarrow (c) \Leftrightarrow (d)$, while $(b)$ does not imply $(a)$.
\end{thm}

We will
prove that $(d)\Rightarrow (c) \Rightarrow (b)
\Leftarrow (a)$ while none of conditions $(a)$, $(b)$ and $(c)$
implies $(d)$ for the weak operator topology (not even for bounded
generators). This immediately shows that also the corresponding
implications in the second Trotter--Kato theorem
\cite[Theorem III.4.9]{engel/nagel:2000} fail for the weak operator
topology. We finish the note with an open question.


\section{The results}

%
%

We first give an example of bounded generators $A$ and $A_n$, $n\in\N$, generating contraction semigroups such that $\lim_{n\to\infty}A_n x=Ax$ weakly for every $x\in X$, but the semigroups $(e^{t A_n})$ do not converge weakly to $(e^{tA})$. This shows that the implication $(a)\Rightarrow (d)$ (as well as $(b)\Rightarrow (d)$) in Theorem \ref{thm:trotter-kato-1} does not hold for the weak operator topology.

\begin{example}\label{ex:counterex-gener}
Consider $X:=l^p$, $1\leq p<\infty$, and the operators $\text{\~A}_n$ defined by
\begin{equation*}
\text{\~A}_n (x_1,\ldots,x_n,x_{n+1},\ldots,x_{2n},\ldots):=(x_{n+1}, x_{n+2},\ldots,x_{2n}, x_1, x_2,\ldots, x_{n},x_{2n+1},x_{2n+2},\ldots)
\end{equation*}
exchanging the first $n$ coordinates of a vector with its next $n$ coordinates.
%
Then $\|\text{\~A}_n\|\leq 1$ implies that $\text{\~A}_n$ generates a $C_0$-semigroup $\text{\~T}_n(\cdot)$ satisfying $\|\text{\~T}_n(t)\|\leq e^t$ for every $n\in \N$ and $t\geq 0$.

The operators $\text{\~A}_n$ converge weakly to zero as $n\to\infty$.
Moreover, $\text{\~A}_n^2=I$ for every $n\in\N$. Therefore
\begin{eqnarray*}
&\text{\~T}_n(t)&=\sum_{k=0}^\infty \frac{t^k \text{\~A}_n^k}{k!}
= \sum_{k=0}^\infty \frac{t^{2k+1}}{(2k+1)!} \text{\~A}_n + \sum_{k=0}^\infty \frac{t^{2k}}{(2k)!} I \\
&\ &= \frac{e^t-e^{-t}}{2}\text{\~A}_n + \frac{e^t+e^{-t}}{2}I
\  \stackrel{\sigma}{\longrightarrow} \frac{e^t+e^{-t}}{2}I \neq I,
\end{eqnarray*}
and this
convergence is uniform on compact time intervals. 
However, the limit does not satisfy the semigroup law.

By rescaling $A_n:=\text{\~A}_n-I$ we obtain a sequence of
contractive semigroups converging weakly  and
uniformly on compact time intervals to a family which is not a
semigroup while the bounded generators converge weakly to $-I$
(which is a generator).
\end{example}

\begin{remark}
The above example shows in particular that (in general) the space of
all contractive $C_0$-semigroups on a Banach
space is not complete for the topology corresponding to the weak
operator convergence, uniform on compact time intervals. Note that
for the strong operator topology this space is complete, a fact
essentially used in Eisner, Ser\'eny \cite{eisner/sereny:2007cont}.
\end{remark}

%

Next we give an example showing that the implication $(c)\Rightarrow (d)$ does not hold for the weak operator topology.
Note that the converse implication is true by the representation $R(\lambda,A_n)=\int_0^\infty e^{-\lambda t}T_n(t)dt$ and  Lebesgue's theorem. Note further that the implications $(a)\Rightarrow (b)$ and $(c)\Rightarrow (b)$ also remain true for the weak operator topology (the first is trivial and the proof of the second is the same as in the strong case).

\begin{example}\label{ex:counterex-res}
On $X:=l^2$ consider the operators
$V_n:=(1-\frac{1}{n})\text{\~A}_n$ for $n\in\N$, where
$\text{\~A}_n$ are the operators from the above example. Then
$\|V_n\|=1-\frac{1}{n}<1$ and hence $1\notin \sigma(V_n)$. By
Foias--Sz.-Nagy \cite[Thm. III.8.1]{sznagy/foias} every $V_n$ is the
cogenerator of a contractive $C_0$-semigroup $S_n(\cdot)$ on $X$.
Recall that the generator $B_n$ of the semigroup $S_n(\cdot)$ is
given by the (negative) Cayley transform of $V_n$, i.e.,
$B_n=(I+V_n)(I-V_n)^{-1}$. In our case, all $B_n$ are bounded
operators as well.

We see that
$R(1,B_n)=\frac{1-V_n}{2}\stackrel{\sigma}{\longrightarrow}
\frac{I}{2} = R(1, -I)$ as $n\to\infty$. Assume that the semigroups
$S_n(\cdot)$ converge weakly and uniformly on compact intervals to
the semigroup $(e^{-t}I)$ generated by $-I$. Then
\begin{equation*}
R^2(1,B_n)=-\frac{d}{d\lambda} R(\lambda, B_n)|_{\lambda=1} =\int_0^\infty e^{-t} t S_n(t)dt \stackrel{\sigma}{\longrightarrow} R^2(1,-I) \text{  as } n\to\infty.
\end{equation*}
However, since $\tilde{A}_n^2=I$ for every $n\in\N$ we obtain
\begin{equation*}
R^2(1,B_n)=\frac{1}{4}\left[I-2V_n+\left(1-\frac{1}{n}\right)^2 I\right]\stackrel{\sigma}{\longrightarrow} \frac{I}{2}\neq R^2(1, -I)
\end{equation*}
which is a contradiction.
\end{example}
%

We now summarise the results shown above.

\begin{thm}\label{thm:trotter-kato-weak}
Let $T(\cdot)$ and $T_n(\cdot)$, $n\in\N$, be $C_0$-semigroups on a Banach space $X$ with ge\-ne\-rators $A$ and $A_n$, $n\in\N$, respectively. Assume that for some constants $M\geq 1$, $\omega\in\R$ the estimates $\|T(t)\|$,  $\|T_n(t)\|\leq Me^{\omega t}$ hold for all $t\geq 0$ and $n\in\N$. Let further $D$ be a core for $A$. Consider the following assertions.
\begin{enumerate}[(a)]
\item $D\subset D(A_n)$ for every $n\in\N$ and $\lim_{n\to\infty}A_n x=Ax$ weakly for all $x\in D$.
\item For every $x\in D$ there exists $x_n\in D(A_n)$ such that $\lim_{n\to\infty}x_n=x$ weakly and \\ $\lim_{n\to\infty}A_n x_n=Ax$ weakly.
\item $\lim_{n\to\infty}R(\lambda, A_n)x=R(\lambda, A)x$ weakly for every $x\in X$ and some $\lambda$ with $\Re\lambda>\omega$.
\item $\lim_{n\to\infty}T_n(t)x=T(t)x$ weakly for every $x\in X$, uniformly in $t$ on compact intervals.
\end{enumerate}
Then $(d)\Rightarrow (c) \Rightarrow (b) \Leftarrow (a)$, while none of conditions $(a)$, $(b)$ and $(c)$ imply $(d)$ in general (even if all operators $A_n$ are bounded).
\end{thm}

\begin{remark}
Let $X$ be a Banach space, $T_n(\cdot)$, $n\in\N$, and $T(\cdot)$ be
$C_0$-semigroups with generators $A_n$, $n\in\N$, and $A$,
respectively, satisfying $\|T(t)\|$, $\|T_n(t)\|\leq Me^{\omega t}$ for some
contants $M\geq 1$, $\omega\in\R$ and all $t\geq 0$ and $n\in\N$,
and take $\lambda$ with $\Re\lambda>\omega$.  If
$\lim_{n\to\infty}T_n(t)
=T(t)$ weakly and uniformly on compact
time intervals, then by the representation $R^k(\lambda,
A)=\frac{1}{(k-1)!}\int_0^\infty e^{(\omega-\lambda) t} t^{k-1}
T_n(t)dt$ we have
\begin{equation}\label{eq:konv-res-one-point}
\lim_{n\to\infty}R^k(\lambda, A_n) = R^k(\lambda, A) \text{ weakly} \ \
\text{for all } k\in\N.
\end{equation}
Note that property (\ref{eq:konv-res-one-point}) for one $\lambda$
implies the same property for all $\lambda$ with $\Re\lambda>\omega$
and the convergence in (\ref{eq:konv-res-one-point}) is uniform for
$\lambda$ in compact sets of the corresponding open right halfplane.
This follows from the power series representation for the resolvent,
a standard connectedness
argument and Cauchy's formula for the
derivatives. Further,
by the same reason (\ref{eq:konv-res-one-point}) is equivalent to
\begin{equation}\label{eq:konv-res}
\lim_{n\to\infty}R(\lambda, A_n) = R(\lambda, A) \text{ weakly and uniformly on compact sets in } \{\lambda:\Re\lambda>\omega\}.
\end{equation}
As we saw in Example \ref{ex:counterex-res}, property
(\ref{eq:konv-res}) is not equivalent to 
weak convergence of the resolvents
 for some (fixed) $\lambda$ as is the case for strong operator convergence.
\end{remark}
Therefore the following question appears naturally.
\begin{question}
 \textit{Is property (\ref{eq:konv-res-one-point}) for some/all $\lambda$ equivalent to weak convergence of the corresponding semigroups uniform on compact time intervals?}
\end{question}
\begin{remark}
If all generators $A_n$ are bounded and satisfy $\|A_n\|\leq M$ for some constant
$M\geq 0$ and all $n\in\N$, then the answer on the above question is ``yes''.
This follows immediately from the Dunford representation
$T_n(t)=\frac{1}{2\pi i} \int_{\{|z|=M+1\}} e^{tz} R(z,A_n) dz$ and (\ref{eq:konv-res}).
\end{remark}

%

\parindent0pt

\end{document}